\documentclass[11pt]{article}
\usepackage{mathtools}
\mathtoolsset{centercolon}
\usepackage{microtype}
\usepackage[usenames,dvipsnames,svgnames,table]{xcolor}
\usepackage{amsfonts,amsmath, amssymb,amsthm,amscd}
\usepackage[height=9in,width=6.5in]{geometry}
\usepackage{verbatim}
\usepackage{tikz}
\usetikzlibrary{decorations.markings}
\tikzset{->-/.style={decoration={
			markings,
			mark=at position .6 with {\arrow{>}}},postaction={decorate}}}
\usepackage{tkz-euclide}

\usepackage[margin=10pt,font=small,labelfont=bf]{caption}
\usepackage{latexsym}
\usepackage[pdfauthor={Hu, Lyons, and Tang},bookmarksnumbered,hyperfigures,colorlinks=true,citecolor=BrickRed,linkcolor=BrickRed,urlcolor=BrickRed,pdfstartview=FitH]{hyperref}
\usepackage{mathrsfs}
\usepackage{graphicx}
\usepackage{enumitem}

\usepackage[inline,nolabel]{showlabels}

\newtheorem{theorem}{Theorem}[section]
\newtheorem{corollary}[theorem]{Corollary}
\newtheorem{lemma}[theorem]{Lemma}
\newtheorem{proposition}[theorem]{Proposition}

\newtheorem{definition}[theorem]{Definition}
\newtheorem{remark}[theorem]{Remark}
\numberwithin{equation}{section}

\newcommand{\be}{\begin{equation}}
	\newcommand{\ee}{\end{equation}}
\newcommand\ba{\begin{align}}
	\newcommand\ea{\end{align}}

\newcommand{\Z}{\mathbb{Z}}

\newcommand\st{\,;\;}
\newcommand\mcP{\mathcal{P}}
\newcommand\mcF{\mathcal{F}}
\newcommand\mcL{\mathcal{L}}
\newcommand\In{\mathsf{In}}
\newcommand\Out{\mathsf{Out}}
\newcommand\rmo{{\mathrm{o}}}

\renewcommand{\P}{\mathbb{P}}

\newcommand{\notion}[1]{{\bf  \textit{#1}}}

\newcommand\und[1]{\overline{\overline {#1}}}
\newcommand\dir[1]{\overrightarrow {#1}}

\DeclareMathOperator{\UST}{\mathsf{UST}}

\newcommand\flr[1]{\lfloor #1 \rfloor}

\providecommand{\keywords}[1]
{
	\small	
	\textbf{\textit{Keywords ---}} #1
}

\begin {document}
\title{A Reverse Aldous--Broder Algorithm}
\author{
	Yiping Hu\thanks{Department of Mathematics, University of Washington. Partially supported by the McFarlan Fellowship from University of Washington. Email: \protect\url{huypken@uw.edu}. }
	\and
		Russell Lyons\thanks{Department of Mathematics, Indiana
		University. Partially supported by the National
		Science Foundation under grants DMS-1612363 and DMS-1954086.
		Email: \protect\url{rdlyons@indiana.edu}.}
	\and
	Pengfei Tang\thanks{Department of Mathematics, Indiana
		University. Partially supported by the National
		Science Foundation under grant DMS-1612363 and by ERC starting grant 676970 RANDGEOM".
      Current affiliation: Department of Mathematical Sciences, Tel Aviv University.
		Email: \protect\url{pengfeitang@mail.tau.ac.il}.}
}
\date{\today}
\maketitle



\begin{abstract}
The Aldous--Broder algorithm provides a way of sampling a uniformly random spanning tree for finite connected graphs  using simple random walk. Namely, start a simple random walk on a connected graph and stop at the cover time. The tree formed by all the first-entrance edges has the law of a uniform spanning tree. Here we show that the tree formed by all the last-exit edges also has the law of a uniform spanning tree. This answers a question of Tom Hayes and Cris Moore from 2010. The proof relies on a bijection that is related to the BEST theorem in graph theory. We also give other applications of our results, including new proofs of the reversibility of loop-erased random walk, of the Aldous--Broder algorithm itself, and of Wilson's algorithm.
\end{abstract}

\begin{abstract}
L'algorithme Aldous--Broder fournit un moyen d'échantillonner un arbre couvrant uniformément aléatoire pour des graphes connexes finis en utilisant une marche aléatoire simple.  À savoir, commençons une marche aléatoire simple sur un graphe connexe et arrêtons-nous au moment de la couverture. L'arbre formé par toutes les arêtes de la première entrée a la loi d'un arbre couvrant uniforme. Nous montrons ici que l'arbre formé par toutes les arêtes de dernière sortie a également la loi d'un arbre couvrant uniforme. Cela répond à une question de Tom Hayes et Cris Moore de 2010. La preuve repose sur une bijection qui est liée au théorème BEST de la théorie des graphes. Nous donnons également d'autres applications de nos résultats, y compris de nouvelles preuves de la réversibilité de la marche aléatoire effacée en boucle, de l'algorithme d'Aldous--Broder lui-même et de l'algorithme de Wilson.
\end{abstract}

\bgroup
\egroup

\section{Introduction}\label{sec:intro}

 Let $G=(V,E)$ be a finite connected graph and let $\mathcal{T}$ be the set of spanning trees of $G$. Obviously $\mathcal{T}$  is a finite set. The \notion{uniform spanning tree} is the uniform measure on $\mathcal{T}$, which is denoted by $\mathsf{UST}(G)$. 
 
 In the late 1980s, Aldous \cite{Aldous1990} and Broder \cite{Broder1989} found an algorithm to generate a uniform spanning tree using simple random walk on $G$; both authors thank Persi Diaconis for discussions. The algorithm is now called the Aldous--Broder algorithm. It  generates a sample of the uniform spanning tree as follows.
 
 Start a simple random walk $( X_n)_{n\geq0}$ on $G$ and stop at the cover time, i.e., stop when it first visits all the vertices of $G$. Collect all the first-entrance edges, i.e., edges $(X_n, X_{n+1})$ for $n \ge 0$ such that $X_{n+1} \ne X_k$ for all $k \le n$. These edges form a random spanning tree, $T$. Then this random tree $T$ has the law $\mathsf{UST}(G)$.  The simple random walk can start from  any vertex of $G$, and so $X_0$ can have any initial distribution.

The proof of this result depends on ideas related to the Markov chain tree theorem. Namely, let $(X_n)_{n \in \Z}$ be a biinfinite stationary random walk on $G$. Let the last exit from a vertex $x$ before time $0$ be $\lambda(x) := \max \{ n < 0 \st X_n = x \}$. Then the random spanning tree $\bigl\{ (x, X_{\lambda(x) + 1}) \st x \in V \setminus \{ X_0\} \bigr\}$ of last-exit edges has the uniform distribution. One then proves the validity of the Aldous--Broder algorithm by reversing time and using reversibility of simple random walk; see, e.g., \cite[Section 4.4]{LP2016} for a proof. The algorithmic aspects were studied by \cite{Broder1989}, while theoretical implications were studied by \cite{Aldous1990} and, later, Pemantle \cite{Pemantle1991}.

 In 2010 (``domino" forum email discussion, 9 Sep.; published in \cite[p.~645]{MooreMertens2011}), Tom Hayes and Cris Moore asked whether the tree formed by  all the last-exit edges starting at time $0$ and ending at the cover time has the law of the uniform spanning tree. Here we give a positive answer to this question.
 
 \begin{theorem}\label{thm: backward AB algorithm}
 	Let $G$ be a finite, connected graph. Start a simple random walk on $G$ and stop at the cover time. Let $T$ be the random tree formed by all the last-exit edges. Then $T$ has the law $\mathsf{UST}(G)$. 
 \end{theorem} 

\begin{remark}
	As we will show, the corresponding result holds for finite, connected networks and the associated network random walk.
\end{remark}

We call this algorithm the \notion{reverse Aldous--Broder} algorithm.

Our proof shows a remarkable strengthening of this equality: Notice that the Aldous--Broder algorithm run up to time $n$ gives a tree $T_n$ on the vertices $\{X_k \st k \le n\}$. The evolution of the rooted tree $(T_n, X_n)$ to $(T_{n+1}, X_{n+1})$ is given by a Markov chain. For this chain, $T_{n+1}$ contains $T_n$ and may be equal to $T_n$. Furthermore, if $\tau_{\textnormal{cov}}$ denotes the cover time, then $T_n = T_{\tau_{\textnormal{cov}}}$ for all $n \ge {\tau_{\textnormal{cov}}}$. Similarly, the reverse Aldous--Broder algorithm run up to time $n$ gives a tree $\overline T_n$ on the vertices $\{X_k \st k \le n\}$. The evolution of $(\overline T_n, X_n)$ to $(\overline T_{n+1}, X_{n+1})$ is also given by a Markov chain. For this chain, $\overline T_{n+1}$ \emph{need not} contain $\overline T_n$. In addition, $\overline T_n$ need not equal $\overline T_{\tau_{\textnormal{cov}}}$ for $n > {\tau_{\textnormal{cov}}}$. Our theorem above is that the distributions of $T_{\tau_{\textnormal{cov}}}$ and $\overline T_{\tau_{\textnormal{cov}}}$ are the same. Our strengthening will show that for every $n$, the distributions of $T_n$ and $\overline T_n$ are the same, even though the two Markov chains have different transition probabilities. Moreover, for every $n$, the distributions of $T_{n \wedge {\tau_{\textnormal{cov}}}}$ and $\overline T_{n \wedge {\tau_{\textnormal{cov}}}}$ are the same. These are both proved as Corollary \ref{cor: same distribution for stopped random walks}. It is also true that that $T_\tau$ has the same distribution as $\overline{T}_\tau$ for some other stopping times $\tau$, such as the hitting time of a vertex; see Proposition \ref{prop: extension to general stopping times}. Using this fact together with the Aldous--Broder algorithm, one obtains the well-known result that the path between two vertices in the $\mathsf{UST}(G)$ has the same distribution as the loop-erased random walk; see Corollary \ref{cor: path in UST has the same law as LERW}.

In outline, our proof that $T_n$ and $\overline T_n$ have the same distribution proceeds as follows. We fix a starting vertex $x$ and an ending vertex $y$ at time $n$. For a path $\gamma$ of length $n$ from $x$ to $y$, there corresponds a first-entrance tree $F(\gamma)$ and a last-exit tree $L(\gamma)$. We define a probability-preserving permutation $\Phi$ of the set of paths of length $n$ from $x$ to $y$ such that $F = L \circ \Phi$. This easily leads to our results. The key, then, is to define $\Phi$. For this purpose, we consider two types of multi-digraphs, one with a linear ordering on each set of incoming edges and the other with a linear ordering on each set of outgoing edges. The former type also determines a spanning tree (of the multi-digraph) directed away from $x$ via first-entrance edges, whereas the latter type determines a spanning tree directed toward $y$ via last-exit edges.  A path $\gamma$ naturally yields a multi-digraph $\hat F(\gamma)$ of the first type and a multi-digraph $\hat L(\gamma)$ of the second type. Reversing all edges of one type of multi-digraph yields one of the other type, except that the edges on the path from $x$ to $y$ in the spanning tree are not reversed. The correspondence of multi-digraphs to paths then allows us to define $\Phi$.

Our proof is closely related to the BEST theorem in graph theory, named after the initials of de Bruijn and van Aardenne-Ehrenfest \cite{BE1951} and Smith and Tutte \cite{ST1941}; we draw the connection below in Theorem \ref{thm: the BEST theorem}. Interestingly, by using the BEST theorem, we are able to provide a new proof of the validity of the Aldous--Broder algorithm. This is done in Theorem \ref{thm: Aldous--Broder}: In outline, we note that the BEST theorem shows that given the transition counts of a random walk path of length $n$, each directed tree is equally likely to correspond to the path. A simple formula connects the transition counts with the likelihood of a tree, and letting $n \to\infty$ yields a spanning tree whose probability is then easily deduced. 

The BEST theorem has been used before to study statistics for Markov chains (see the survey by Billingsley \cite{Billingsley1961}) and to characterize Markov exchangeable sequences, as observed by Zaman \cite{Zaman1984}.

\section{Proof of the Main Theorem}\label{sec: 2}

We allow weights $w\colon E \to (0, \infty)$ on the edges, so that $(G, w)$ is a network, and the corresponding \notion{weighted uniform spanning tree} measure puts mass proportional to $\prod_{e \in T} w(e)$ on a spanning tree $T$. We could also allow parallel edges and loops in $G$, but this would simply require more complicated notation, so we assume that $G$ is a simple graph. Write $w(x)$ for the sum of $w(e)$ over all edges $e$ incident to $x$.

Let $x$ and $y$ be two vertices of $G$. If there is an edge in $G$ connecting $x$ and $y$, then we write $x\sim y$. We call $\gamma=(v_0,\ldots,v_n)$  a \notion{path} (or a walk) on $G$ from $v_0$ to $v_n$ if $v_{i-1}\sim v_i$ for $i=1,\ldots,n$, and $|\gamma|:=n$ is called the \notion{length} of the path $\gamma$.

Let  $\mathcal{P}_n^{x,y}$ denote the set of all paths in $G$ from $x$ to $y$ with length $n$. Simple random walk is the Markov chain on $V$ with transition probabilities $p(x, y) :=\frac{1}{\textnormal{deg}(x)}\mathbf{1}_{\{x\sim y \}}$. The network random walk on $(G, w)$ has instead the transition probabilities $p(x, y) := \frac{w(x, y)}{w(x)}\mathbf{1}_{\{x\sim y \}}$. All results quoted earlier hold for network random walks, not only simple random walks, and for the corresponding random spanning trees. In particular, the Aldous--Broder algorithm holds in this more general context: the network random walk generates the weighted uniform spanning tree by collecting the first-entrance edges up to the cover time \cite[Section 4.4]{LP2016}.
We henceforth assume that our graphs are weighted and use the corresponding network random walks; the weighted uniform spanning tree measure is still denoted by $\UST(G)$.

Let $\mathscr{T}$ be the set of all subtrees of $G$, including those that are not necessarily spanning. 

For a path $\gamma=(v_0,\ldots,v_n)\in \mathcal{P}_n^{x,y}$, where $v_0 = x$ and $v_n = y$, we write $V(\gamma)=\{v_0, \ldots, v_n\}$ for the set of vertices of $\gamma$.  
For each $u\in V(\gamma)\backslash\{ x\}$, there is a smallest index $i\geq 1$ such that $v_i=u$; we call the edge $(v_{i-1},v_i)$ the \notion{first-entrance edge} to $u$.  

Define the \notion{first-entrance operator} $F\colon\mathcal{P}_n^{x,y}\rightarrow \mathscr{T}$ by setting $F(\gamma)$ to be the tree formed by all the first-entrance edges to vertices in $V(\gamma)\backslash\{ x\}$.

Similarly, for each $u\in V(\gamma)\backslash\{ y\}$, there is a largest index $i\leq n-1$ such that $v_i=u$, and we call the edge $e=(v_{i},v_{i+1})$ the \notion{last-exit edge} of $u$. 
Define the \notion{last-exit operator} $L\colon\mathcal{P}_n^{x,y}\rightarrow \mathscr{T}$ by setting $L(\gamma)$ to be the tree formed by all the last-exit edges of vertices in $V(\gamma)\backslash\{ y\}$.

For a finite path $\gamma = (v_0, \ldots, v_n)$ in $G$ with length $n$ started at $v_0$, we write 
\begin{equation} \label{eq:wts}
p(\gamma):=\mathbb{P}_{v_0}\bigl[(X_0,\ldots,X_n)=\gamma\bigr]
=
\frac{\prod_{k=0}^{n-1} w(v_k, v_{k+1})}{\prod_{k=0}^{n-1} w(v_k)},
\end{equation}
where $\mathbb{P}_{v_0}$ denotes the law of the network random walk on $(G, w)$ started from $v_0$.

The proof will go as follows: In Section \ref{sec: Main Lemma}, we state the key lemma (Lemma \ref{lem: bijection}) and use it to derive Theorem \ref{thm: backward AB algorithm}, together with other consequences. In Section \ref{sec: operators}, we start the proof of the key lemma by establishing the correspondence between random walk paths and colored multi-digraphs. We finish the proof in Section \ref{sec: proof of Main Lemma}.

\subsection{The Main Lemma and its Consequences}
\label{sec: Main Lemma}

The main ingredient for proving our theorem is the following lemma, which can be viewed as an extension of Proposition 2.1 of \cite{Lawler1983}.

\begin{lemma}\label{lem: bijection}
	For every $x,y\in V(G)$ and $n\geq 0$, there is a bijection $\Phi\colon\mathcal{P}_n^{x,y}\rightarrow\mathcal{P}_n^{x,y}$ such that $F=L\circ \Phi$ and 
	\be\label{eq: preserve the random walk measure}
	\forall \gamma \in \mathcal{P}_n^{x,y} \qquad p(\gamma)=p\bigl(\Phi(\gamma)\bigr).
	\ee
   Moreover, $\Phi$ preserves the number of times each vertex is visited as well as the number of times that each (unoriented) edge is crossed.
\end{lemma}

It follows that $\Phi$ is also a measure-preserving bijection on each subset of $\mcP_n^{x, y}$ specified by how many times each edge is crossed or which vertices are visited. For example, such a subset is the set of paths such that every vertex is visited and $y$ is visited only once and only after all other vertices are visited.

Before proving Lemma \ref{lem: bijection}, we show how it gives our claims.

Consider the network random walk $(X_n)_{n=0}^\infty$ on $G$ with arbitrary initial distribution, and recall that $\tau_{\textnormal{cov}}$ denotes the cover time of the graph $G$. Write $U \overset{\mathscr{D}}{=} V$ when $U$ and $V$ have the same distribution. 

\begin{corollary}\label{cor: same distribution for stopped random walks}
	We have for all $n \in \mathbb{N}$,
	\begin{enumerate}[label=(\roman*)]
		\item $F\bigl((X_0,\ldots,X_n)\bigr) \overset{\mathscr{D}}{=} L\bigl((X_0,\ldots,X_n)\bigr)$ and
		\item $F\bigl((X_0,\ldots,X_{\tau_{\textnormal{cov}}\wedge n})\bigr) \overset{\mathscr{D}}{=} L\bigl((X_0,\ldots,X_{\tau_{\textnormal{cov}}\wedge n})\bigr)$. 
	\end{enumerate}
\end{corollary}

\begin{proof}
This is immediate from Lemma \ref{lem: bijection} and the remarks following it.
\end{proof}

\begin{proof}[Proof of Theorem \ref{thm: backward AB algorithm}]
	Letting $n\rightarrow\infty$ in Corollary \ref{cor: same distribution for stopped random walks}(ii)
	and noting that $\mathbb{P}[\tau_{\textnormal{cov}}\geq n]\rightarrow0 $ as $n\rightarrow\infty$, one obtains that 
	\be\label{eq: same distribution at cover time}
	F\bigl((X_0,\ldots,X_{\tau_{\textnormal{cov}}})\bigr) \overset{\mathscr{D}}{=} L\bigl((X_0,\ldots,X_{\tau_{\textnormal{cov}}})\bigr).
	\ee  
	By the Aldous--Broder algorithm, $F\bigl((X_0,\ldots,X_{\tau_{\textnormal{cov}}})\bigr)$ has the law of $\mathsf{UST}(G)$, and thus so does $L\bigl((X_0,\ldots,X_{\tau_{\textnormal{cov}}})\bigr)$.
\end{proof}

	The equation \eqref{eq: same distribution at cover time} holds for more general stopping times (but not all stopping times). Here we give a sufficient condition.

	 Let $C_e^n:=\#\bigl\{j\st 0\leq j\leq n-1,\, \{X_j,X_{j+1}\}=e\bigr\}$ denote the total number of crossings of the edge $e$ in both directions up to time $n$. If $I_y := \sum_{e \ni y} C_e^n$, then the number of visits to $y$ by time $n$ is $\flr{(I_y+1)/2}$. Furthermore, $y \ne x$ is visited at time $n$ iff $I_y$ is odd, and $x$ is visited at time $n$ iff $I_x$ is even. If $\tau$ is a finite stopping time such that 
    \be\label{occup-stop}
    \forall n\geq 0\quad \{\tau=n\}\in\sigma\big(C_e^n\st  e\in E(G)\big),
    \ee
      then $\tau$ satisfies the corresponding version of equation  \eqref{eq: same distribution at cover time}:
	\be\label{eq: general stopping time}
	F\bigl((X_0,\ldots,X_{\tau})\bigr) \overset{\mathscr{D}}{=} L\bigl((X_0,\ldots,X_{\tau})\bigr).
	\ee
For example, if $K$ is  a finite subset of $V(G)$ and $\tau$ is the hitting time of $K$, then $\tau$ satisfies \eqref{eq: general stopping time}.

This implies that the path between two fixed vertices in the uniform spanning tree of a graph has the same law as the loop-erased random walk from one of those vertices to the other. We state this as
Corollary \ref{cor: path in UST has the same law as LERW} below, which was first proved by Pemantle \cite{Pemantle1991} using Proposition 2.1 of Lawler  \cite{Lawler1983}. 	Using Corollary \ref{cor: path in UST has the same law as LERW} and the spatial Markov property of $\mathsf{UST}(G)$, one can immediately get Wilson's algorithm \cite{Wilson1996} for sampling $\mathsf{UST}(G)$ (in networks). This was observed by Wilson \cite{Wilson1996} and Propp and Wilson \cite{ProppWilson1998}. 
We first recall the definition of loop-erased random walk. 
\begin{definition}\label{def: loop-erasure of a finite path}
	Let $\gamma=(v_0,v_1,\ldots,v_n)$ be a finite path on $G$. We define the \notion{loop-erasure} of $\gamma$ to be the self-avoiding path $\textnormal{LE}(\gamma)=(u_0,\ldots,u_m)$ obtained by erasing cycles on $\gamma$ in the order they appear. Equivalently, we set $u_0:=v_0$ and $l_0:=\max\{j\st  0\leq j\leq n,v_j=u_0 \}$. If $l_0=n$, we set $m:=0$, $\textnormal{LE}(\gamma):=(u_0)$ and terminate. Otherwise, we set $u_1:=v_{l_0+1}$ and $l_1:=\max\{j\st   l_0\leq j\leq n, v_j=u_1\}$. If $l_1=n$, we set $m:=1$, $\textnormal{LE}(\gamma):=(u_0,u_1)$ and terminate. Continue in this way: while $l_k<n$, set $u_{k+1}:=v_{l_k+1}$ and $l_{k+1}:=\max\{j\st   l_k\leq j\leq n, v_j=u_{k+1}\}$. Since $\gamma$ is a finite path, we must terminate after at most $n$ steps and get the loop-erased path $\textnormal{LE}(\gamma)$.
\end{definition}

\begin{definition}\label{def: loop-erased random walk}
	Suppose that $u$ and $v$ are distinct vertices in the graph $G$. Start a network random walk on $G $ from $u$ and stop at the first visit of $v$. The loop-erasure of this random walk path is called the \notion{loop-erased random walk} from $u$ to $v$. 
\end{definition}

\begin{corollary}\label{cor: path in UST has the same law as LERW}
	Let $x$ and $y$ be two distinct vertices in $G$, and let $\gamma_{x,y}$ be the unique path in the $\mathsf{UST}(G)$ from $x$ to $y$. Then $\gamma_{x,y}$ has the same law as the loop-erased random walk from $x$ to $y$. In particular, the loop-erased random walk from $x$ to $y$ has the same law as the loop-erased random walk from $y$ to $x$.
\end{corollary}
\begin{proof}
	Start a random walk on $G$ from $x$. The stopping time $\tau:=\min\{k\st  X_k=y\}$ satisfies \eqref{occup-stop}, whence also \eqref{eq: general stopping time}.
	By the Aldous--Broder algorithm, $\gamma_{x,y}$ has the same law as the unique path in $F\bigl((X_0,\ldots,X_{\tau})\bigr)$ from $x$ to $y$. Thus   $\gamma_{x,y}$ has the same law as the unique path in $L\bigl((X_0,\ldots,X_{\tau})\bigr)$ from $x$ to $y$, which is just the loop-erased random walk from $x$ to $y$.

The last sentence follows from the equality $\gamma_{x, y} = \gamma_{y, x}$.
\end{proof}

Stopping times that rely on more randomness can also satisfy \eqref{eq: general stopping time}. One interesting such case is to consider a geometric random variable $T$ independent of the random walk and let $K$ be a finite subset of $V(G)$. Define $\tau$ to be the time of the $T$th visit to  $K$. Then $\tau$ is such a stopping time that satisfies \eqref{eq: general stopping time}. 

\begin{proposition}\label{prop: extension to general stopping times}
	If $T$ is a random time independent of the random walk, and $\tau$ is a finite stopping time such that $\{\tau=n\}\in\sigma\big(T,\, C_e^n\st  e\in E(G)\big)$ for all $n\geq 0$, then $\tau$ satisfies \eqref{eq: general stopping time}.
\end{proposition}
\begin{proof}
	
 For each path $\gamma=(v_0,\ldots,v_n)$ in $G$, we claim that 
\be\label{eq: decomposition into paths}
\mathbb{P}\big[(X_0,\ldots,X_n)=\gamma,\, \tau=n\big]=\mathbb{P}\big[(X_0,\ldots,X_n)=\Phi(\gamma),\, \tau=n\big].
\ee
By Lemma \ref{lem: bijection} and the independence between $T$ and the random walk, we have for every $k \ge 0$,
\be\label{eq: indep between T and SRW}
\mathbb{P}\big[(X_0,\ldots,X_n)=\gamma,\, T=k\big]=\mathbb{P}\big[(X_0,\ldots,X_n)=\Phi(\gamma),\, T=k\big].
\ee
Let $m_e$ denote the number of crossings of the edge $e$ in $\gamma$.
Since $\Phi$ preserves the number of times each edge is crossed, 
\[
\bigl\{(X_0,\ldots,X_n)=\gamma,\, T=k\bigr\} \cup \bigl\{(X_0,\ldots,X_n)=\Phi(\gamma),\, T=k\bigr\}
\subset \{ T=k,\, \forall{e\in E}\enspace C_e^n=m_e\}.
\]
By the assumption on $\tau$, 
\[
\{ T=k,\, \forall{e\in E}\enspace C_e^n=m_e\}  \cap\{\tau=n\}=\varnothing
\quad\text{or}\quad 
\{ T=k,\, \forall{e\in E}\enspace C_e^n=m_e\}\subset\{\tau=n\}.
\]
In the former case,
\be\label{eq: to be summed over k}
\mathbb{P}\big[(X_0,\ldots,X_n)=\gamma,\, T=k,\, \tau=n\big]=\mathbb{P}\big[(X_0,\ldots,X_n)=\Phi(\gamma),\, T=k,\, \tau=n\big]=0.
\ee
In the latter case,
\begin{eqnarray}
\mathbb{P}\big[(X_0,\ldots,X_n)=\gamma,\, T=k,\, \tau=n\big]&=&
\mathbb{P}\big[(X_0,\ldots,X_n)=\gamma,\, T=k\big]\nonumber\\
&\stackrel{\eqref{eq: indep between T and SRW}}{=}&\mathbb{P}\big[(X_0,\ldots,X_n)=\Phi(\gamma),\, T=k\big]\nonumber\\
&=&\mathbb{P}\big[(X_0,\ldots,X_n)=\Phi(\gamma),\, T=k,\, \tau=n\big].
\end{eqnarray}
So we always have
\[
\mathbb{P}\big[(X_0,\ldots,X_n)=\gamma,\, T=k,\, \tau=n\big]=\mathbb{P}\big[(X_0,\ldots,X_n)=\Phi(\gamma),\, T=k,\, \tau=n\big].
\]
Summing this equality over $k$ yields \eqref{eq: decomposition into paths}.

Now it is easy to see that \eqref{eq: general stopping time} holds for such stopping times $\tau$: for each subtree $t$ of $G$, one has that 
\begin{eqnarray*}
\mathbb{P}\big[	F\bigl((X_0,\ldots,X_{\tau})\bigr)=t\big]&=&
\sum_{\gamma:F(\gamma)=t}\mathbb{P}\big[(X_0,\ldots,X_{\tau})=\gamma,\, \tau=|\gamma|\big]\\
&\stackrel{\eqref{eq: decomposition into paths}}{=}&\sum_{\gamma:F(\gamma)=t}\mathbb{P}\big[(X_0,\ldots,X_{\tau})=\Phi(\gamma),\, \tau=|\gamma|\big]\\
&\stackrel{\gamma'=\Phi(\gamma)}{=}&\sum_{\gamma':L(\gamma')=t}\mathbb{P}\big[(X_0,\ldots,X_{\tau})=\gamma',\, \tau=|\gamma'|\big]\\
&=&\mathbb{P}\big[	L\bigl((X_0,\ldots,X_{\tau})\bigr)=t\big]. 
\end{eqnarray*}
      \par          \vspace{-1.6\baselineskip} 
\end{proof}
\vspace{1\baselineskip}

\subsection{Augmented Operators \texorpdfstring{$\hat L$ and $\hat F$}{} and the Reversing Operator $R$}
\label{sec: operators}

Note that if $x = y$, then Lemma \ref{lem: bijection} is trivial: take $\Phi$ to reverse the path. Thus, we assume from now on that $x \ne y$.

We encode every path by a colored multi-digraph. We say a multi-digraph marked with a \notion{start} vertex $x$ and an \notion{end} vertex $y$ is \notion{balanced} if the indegree of each vertex $u \notin \{x, y\}$ equals the outdegree of $u$,
the outdegree of $x$ is larger than its indegree by one, and the indegree of $y$ is larger than its outdegree by one.
The balance property is necessary for the existence of an Eulerian path in a multi-digraph.
The coloring is either of the following two types.

By an \notion{exit coloring} of a multi-digraph, we mean an assignment to each vertex $v$ of a linear ordering on the set $\Out(v)$ of its outgoing edges. We call an edge \notion{lighter} or \notion{darker} if it is smaller or larger in the ordering. The maximal (darkest) edge in $\Out(v)$ is regarded as \notion{black}, except for $v = y$. For the \notion{entrance coloring}, the ordering is on the edges $\In(v)$ incoming to $v$ instead, and 
there is no black edge leading into $x$.

Let $\mathcal{L}_n^{x,y}$ (resp., $\mathcal{F}_n^{x,y}$) denote the set of all colored multi-digraphs marked with the start $x$ and end $y$ satisfying the following five properties:
\begin{enumerate}[label=(\roman*)]
	\item the vertex set of the multi-graph contains the start $x$ and end $y$ and is a subset of the vertex set of the original graph $G$;
	\item there are $n$ directed edges in total, each one being an orientation of an edge from $G$;
	\item the multi-digraph is balanced;
	\item the coloring is an exit coloring (resp., an entrance coloring);
	\item all black edges form a directed spanning tree of the multi-digraph with all edges going from leaves to the root $y$ (resp., from the root $x$ to leaves).
   \label{SpanTreeCond}
\end{enumerate}

%

We define the augmented last-exit operator $\hat L\colon \mathcal{P}_n^{x,y} \to \mathcal{L}_n^{x,y}$ that, given a path, returns a colored multi-digraph in the most obvious way: the vertex set consists of all the visited vertices, and we draw one directed edge for each step of the path. This multi-digraph is clearly balanced. The exit coloring is naturally given by the order in which these outgoing edges are traversed. Also recall that there is no black edge leading out of $y$. Note that the tree formed by all black edges, ignoring their orientation, is exactly the output of the last-exit operator $L$ defined above, so $\hat L$ can be viewed as an augmentation of $L$. All conditions can be checked easily, so the map is well defined. 

Similarly, we can define the augmented first-entrance operator $\hat F\colon \mathcal{P}_n^{x,y} \to \mathcal{F}_n^{x,y}$. The only difference is that we consider the entrance coloring on incoming edges instead of the exit coloring on outgoing edges, and we use the \emph{reverse} of the order in which incoming edges are traversed, so that black edges, in particular, are first-entrance edges. 

Clearly $\hat L$ and $\hat F$ are both injections.	

\begin{lemma}
\label{lem: BEST bijection}
	The maps $\hat L$ and $\hat F$ are bijections.	
\end{lemma}

\begin{proof}
	We describe the inverse operator for $\hat L$. Given a colored multi-digraph in $\mathcal{L}_n^{x,y}$, we can associate it with a path by starting from $x$ and traversing the multi-digraph according to the following instruction: $$\text{always exit from the lightest unused edge until you run out of edges to use.}$$
	Notice that by the balance property, the algorithm always terminates at $y$. It remains to show that the path is Eulerian (i.e., uses all $n$ edges exactly once), so that the operator is well defined. Given this, since the above rule respects the exit coloring, we may deduce that this traversal algorithm gives the inverse, $\hat L^{-1}$.

	If the path is not Eulerian, i.e., if some edge has not been used, then some outgoing edge from some vertex $u$ has not been used. According to our
   instructions, it follows that the black edge leaving $u$ has not been used. Let that black edge be $(u, v)$. By the balance condition, some outgoing edge from
   $v$ has not been used, whence the black edge leaving $v$ has not been used (unless $v = y$). By repeating this argument and using condition \ref{SpanTreeCond}, we arrive at the conclusion that some edge leaving $y$ has not been used, whence the path cannot have terminated, a contradiction.

   The proof for $\hat F$ is similar; it also follows from \eqref{defPhi} below and the fact that the other functions there are bijections.
\end{proof}

The bijection $\hat L$ was also used in the proof of the BEST theorem; see Theorem \ref{thm: the BEST theorem} and its proof for details, where we use $\hat F$ instead.

	Define a reversing operator $R\colon \mathcal{F}_n^{x,y} \to \mathcal{L}_n^{x,y}$ that maintains the vertex set but reverses all edges, except that the directions of the black edges on the unique path from $x$ to $y$ in the spanning tree remain unchanged. The coloring is also maintained: for $H \in \mcF_n^{x, y}$ and a vertex $u$, the set $\In(u)$ in $H$ is mapped to $\Out(u)$ in $R(H)$, and so this set of edges can maintain its linear ordering, except that if $u \notin \{x, y\}$ is on the path from $x$ to $y$, then one black edge is replaced by another, whereas $x$ gains an outgoing black edge in $R(H)$ compared to the incoming edges to $x$ in $H$, and $y$ loses an outgoing black edge in $R(H)$ compared to the incoming edges to $y$ in $H$. It is not hard to see that the operator $R$ does have codomain $\mcL_n^{x, y}$ and is bijective. 

\subsection{Proof of Lemma \ref{lem: bijection} }
\label{sec: proof of Main Lemma}

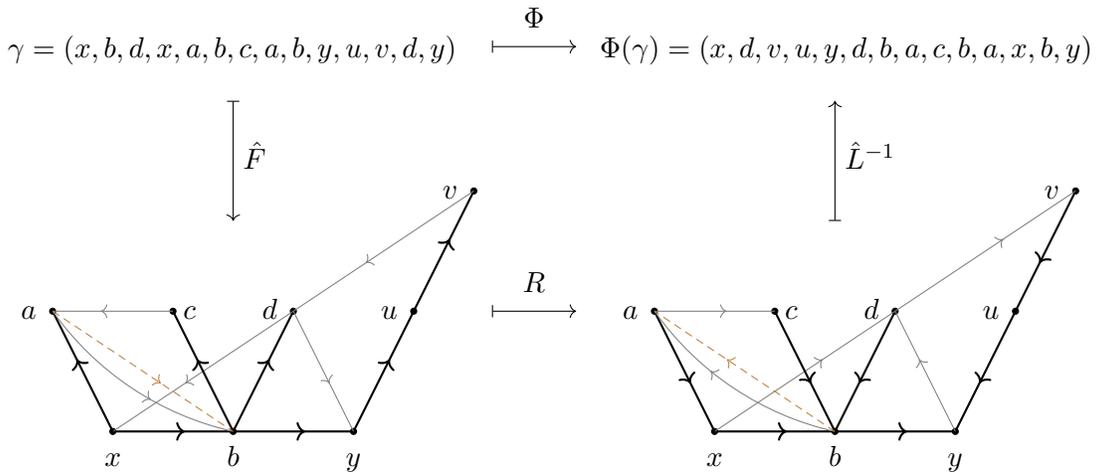
\begin{figure}[h!]
      \centering
		\begin{tikzpicture}[scale=0.8, text height=1.5ex,text depth=.25ex] 

		\draw[fill=black] (1,0) circle [radius=0.05]; 
		\node[below] at (1,-0.1) {$x$};
		
		\draw[fill=black] (0,2) circle [radius=0.05]; 
		\node[left] at (-0.1,2) {$a$};

		\draw[fill=black] (3,0) circle [radius=0.05]; 
		\node[below] at (3,-0.1) {$b$};
		
		\draw[fill=black] (5,0) circle [radius=0.05]; 
		\node[below] at (5,-0.1) {$y$};

		\draw[fill=black] (2,2) circle [radius=0.05]; 
		\node[right] at (2,2) {$c$};
		
		\draw[fill=black] (4,2) circle [radius=0.05]; 
		\node[left] at (4-0.1,2) {$d$};
		
		\draw[fill=black] (6,2) circle [radius=0.05]; 
		\node[left] at (6-0.1,2) {$u$};
		
		\draw[fill=black] (7,4) circle [radius=0.05]; 
		\node[left] at (7-0.1,4) {$v$};
		
		\draw [->-,thick] (1,0)--(3,0); 
		\draw [->-,thick] (3,0)--(4,2); 
		\draw [->-,color=gray] (4,2)--(1,0); 
		\draw [->-,thick] (1,0)--(0,2); 
         \begin{scope}[decoration={markings,
			   mark=at position .6 with {\arrow{>}}}]
        \tkzDefPoint(0,2){a}\tkzDefPoint(3,0){b}\tkzDefPoint(1.2,.8){C}
        \tkzCircumCenter(a,b,C)\tkzGetPoint{O}
        \tkzDrawArc[postaction={decorate},color=gray](O,a)(b)
         \end{scope}

		\draw [->-,thick] (3,0)--(2,2); 
		\draw [->-,color=gray] (2,2)--(0,2); 
		
		\draw [->-,color=brown,densely dashed] (0,2)--(3,0); 
		
		\draw [->-,thick] (3,0)--(5,0); 
		\draw [->-,thick] (5,0)--(6,2); 
		\draw [->-,thick] (6,2)--(7,4); 
		\draw [->-,color=gray] (7,4)--(4,2); 
		\draw [->-,color=gray] (4,2)--(5,0); 

		
		\draw[fill=black] (1+10,0) circle [radius=0.05]; 
		\node[below] at (1+10,-0.1) {$x$};
		
		\draw[fill=black] (0+10,2) circle [radius=0.05]; 
		\node[left] at (-0.1+10,2) {$a$};

		\draw[fill=black] (3+10,0) circle [radius=0.05]; 
		\node[below] at (3+10,-0.1) {$b$};
		
		\draw[fill=black] (5+10,0) circle [radius=0.05]; 
		\node[below] at (5+10,-0.1) {$y$};

		\draw[fill=black] (2+10,2) circle [radius=0.05]; 
		\node[right] at (2+10,2) {$c$};
		
		\draw[fill=black] (4+10,2) circle [radius=0.05]; 
		\node[left] at (4+10-0.1,2) {$d$};
		
		\draw[fill=black] (6+10,2) circle [radius=0.05]; 
		\node[left] at (6+10-0.1,2) {$u$};
		
		\draw[fill=black] (7+10,4) circle [radius=0.05]; 
		\node[left] at (7+10-0.1,4) {$v$};

		\draw [->-,thick] (1+10,0)--(3+10,0); 
		\draw [->-,thick] (4+10,2)--(3+10,0); 
		\draw [->-,color=gray] (1+10,0)--(4+10,2); 
		\draw [->-,thick]   (10,2)--(1+10,0);  
         \begin{scope}[decoration={markings,
			   mark=at position .4 with {\arrow{<}}}]
        \tkzDefPoint(10,2){a}\tkzDefPoint(13,0){b}\tkzDefPoint(11.2,0.8){C}
        \tkzCircumCenter(a,b,C)\tkzGetPoint{O}
        \tkzDrawArc[postaction={decorate},color=gray](O,a)(b)
         \end{scope}
		\draw [->-,thick]  (2+10,2)--(3+10,0); 
		\draw [->-,color=gray]  (10,2)--(2+10,2); 
		
		\draw [->-,color=brown,densely dashed] (3+10,0)--(10,2);  
		
		\draw [->-,thick] (3+10,0)--(5+10,0);  
		\draw [->-,thick] (6+10,2)--(5+10,0);  
		\draw [->-,thick]  (7+10,4)--(6+10,2);  
		\draw [->-,color=gray] (4+10,2)--(7+10,4);   
		\draw [->-,color=gray]   (5+10,0)--(4+10,2); 

		
		\draw  [|->] (7.3,2)--(8.7,2);
		\node [above] at (8,2.1) {$R$};

		\draw  [|->] (7.3,6.4)--(8.7,6.4);
		\node [above] at (8,6.5) {$\Phi$};

		\node[above] at (3,6) {$\gamma=(x,b,d,x,a,b,c,a,b,y,u,v,d,y)$};
		
		\draw  [|->] (3,5.5)--(3,3.5);
		\node [right] at (3,4.5) {$\hat F$};

		\node[above] at (3+10.2,6) {$\Phi(\gamma)=(x,d,v,u,y,d,b,a,c,b,a,x,b,y)$};
		\draw  [|->] (3+10,3.5)--(3+10,5.5);
		\node [right] at (3+10,4.5) {$\hat L^{-1}$};
		
		\end{tikzpicture}
		
		\caption{An example of the bijection $R$, where we use the order that brown, dashed edges are lighter than gray edges and gray edges are lighter than black edges. In addition, black edges are drawn thick; these form the spanning tree. Note that $R$ reverses the orientation of every edge except for the black edges that form a path from $x$ to $y$.}
	\label{fig: bijection R}  
   \end{figure}

\begin{proof}[Proof of Lemma \ref{lem: bijection}]
	We may now define the main construction, illustrated in Figure \ref{fig: bijection R}:
   \begin{equation}  \label{defPhi}
   \Phi := \hat{L}^{-1} \circ R \circ \hat{F}.
   \end{equation}
	Since $R$ maintains the coloring, the black edges that formed the first-entrance spanning tree in $\mathcal{F}_n^{x,y}$ now form the last-exit version in $\mathcal{L}_n^{x,y}$ after reversal, that is, $F = L \circ \Phi$. Because $\hat F$, $R$, and $\hat L^{-1}$ are injections, so is $\Phi$, which forces $\Phi$ to be a bijection. The map $\Phi$ preserves the random walk measure because all edges and vertices are visited same number of times: see \eqref{eq:wts}.
 \end{proof}

\section{Concluding Remarks}

\subsection{The BEST Theorem and the Aldous--Broder Algorithm}

We will make explicit the relationship of our proof to the BEST theorem, and use the latter to give an alternative proof of the validity of the Aldous--Broder algorithm.
In addition, because of their possible interest, we show how similar ideas lead to new proofs of the Markov chain tree theorem and the directed version of Wilson's algorithm, although these are not the shortest known proofs.
For additional proofs that are somewhat related to ours and history of many of these results, see \cite{PitmanTang}.

Suppose $D$ is a multi-digraph marked with a start  vertex $x$ and an end vertex $y$. An \notion{Eulerian path} in $D$ is a directed path that starts at $x$, ends at $y$, and uses each edge exactly once. An \notion{arborescence} in $D$ is a directed spanning tree of $D$ with each edge belonging to a directed path from the root $x$ to a leaf. Each Eulerian path gives rise to the arborescence formed by all the first-entrance edges of vertices in $D$ except for that of $x$. 

The original BEST theorem \cite[Theorem 5b]{BE1951} was concerned with \notion{Eulerian circuits}, which are Eulerian paths that start and end at the same vertex. Depending on the context, they are viewed as loops either with or without a distinguished start or end. Here we state a version of the theorem that holds more generally for Eulerian paths. We denote the indegree of a vertex $v$ by $\textnormal{indeg}(v)$.

\begin{theorem}[BEST]\label{thm: the BEST theorem}
	Suppose $D$ is a balanced multi-digraph marked with a start $x$ and an end $y$. For each arborescence $\Lambda$ in $D$, there are exactly
	\be \label{eq: number of Eulerian paths}
		\textnormal{indeg}(x)\prod_{v \in V} (\textnormal{indeg}(v)-1)!
	\ee
	Eulerian paths in $D$ whose first-entrance edges coincide with the arborescence $\Lambda$. In particular, the total number of Eulerian paths in $D$ is given by
	$$
		t_x(D)\,\textnormal{indeg}(x)\prod_{v \in V} (\textnormal{indeg}(v)-1)!,
	$$
	where $t_x(D)$ is the number of arborescences in $D$ directed away from the root $x$.
\end{theorem}

\begin{proof}
	It suffices to show the formula in (\ref{eq: number of Eulerian paths}). The proof is a direct consequence of the bijection $\hat F$ (see Lemma \ref{lem: BEST bijection}): the number of Eulerian paths with first-entrance edges being $\Lambda$ is equal to the number of entrance colorings of $D$ whose black edges correspond to $\Lambda$, while the latter is exactly (\ref{eq: number of Eulerian paths}).
\end{proof}

\begin{remark}
	It is perhaps more natural to reverse the Eulerian paths and formulate Theorem \ref{thm: the BEST theorem} in terms of last-exit edges. However, the first-entrance version will be useful for us to study the Aldous--Broder algorithm in Theorem \ref{thm: Aldous--Broder}.
\end{remark}

\begin{theorem}[Aldous--Broder]
\label{thm: Aldous--Broder}
	Let $G$ be a finite, connected network. Start a network random walk $(X_n)_{n=0}^\infty$ on $G$ from the vertex $x$ and stop at the cover time. Let $T$ be the random spanning tree formed by all the first-entrance edges. Then $T$ has the law $\mathsf{UST}(G)$. 
\end{theorem}

\begin{proof}
Recall that the Aldous--Broder algorithm run up to time $n$ gives a tree on the vertices $\{X_k \st k \le n\}$. We denote this tree by $T_n$. It suffices to show that for every spanning tree $t$,
	\be
	\label{eq: spanning tree probability}
	\lim_{n\to\infty} \mathbb{P}_x\bigl[T_n=t\bigr] =\prod_{e \in t} w(e) \bigg/ \sum_{t'} \prod_{e \in t'} w(e), 
	\ee
	where the summation is over all spanning trees of $G$.
	
	For any pair of vertices $(u, v)$, let
	\be
	\label{eq: transition counts}
	N_{u, v}(n) :=\#\{j\st 0\leq j\leq n-1,\, X_j=u,\, X_{j+1}=v\}
	\ee
	denote the transition counts from $u$ to $v$ up to time $n$. Let $N(n) \colon (u, v) \mapsto N_{u, v}(n)$ be the collection of such transition counts. Note that $N(n)$ uniquely determines the end vertex $X_n$. Based on $N(n)$, we define a balanced multi-digraph $D_n$ with start $X_0=x$ and end $X_n$ that has the same vertex set as $G$ and contains exactly $N_{u, v}(n)$ directed edges from $u$ to $v$ for all pairs $(u, v)$. By definition, all random walk paths of length $n$ with the same transition counts $N(n)$ occur with equal probability. Therefore, the random walk path up to time $n$, conditioned on the transition counts $N(n)$, has the same distribution as the law $\mathbb{Q}_{D_n}$ of a uniformly random Eulerian path in the multi-digraph $D_n$. Let $\dir G$ be the directed graph obtained from $G$ by replacing each edge with two oppositely oriented edges. For any subdigraph $H$ of $D_n$, define $\und H$ to be the subdigraph of $\dir G$ formed from $H$ by identifying each edge of $D_n$ with its corresponding edge in $\dir G$. Thus we have
	\be
	\label{eq: conditional law on transition counts}
		\mathbb{P}_x\bigl[T_n=t \mid N(n)\bigr] = \mathbb{Q}_{D_n}\Bigl[ \und{\Lambda_n} = \Lambda(t)\Bigr], 
	\ee
	where $\Lambda_n$ represents the arborescence generated by the random Eulerian path, and $\Lambda(t)$ is the directed spanning tree in $\dir G$ if we view $t$ as directed away from the root $x$.
	
	By Theorem \ref{thm: the BEST theorem}, each arborescence is equally likely to occur under $\mathbb{Q}_{D_n}$. Moreover, given a spanning tree $t$, there are exactly $\prod_{(u, v) \in \Lambda(t)} N_{u, v}(n)$ arborescences $\Lambda$ in $D_n$ with $\und \Lambda = \Lambda(t)$. Thus, the right-hand side of (\ref{eq: conditional law on transition counts}) is given by
	\be
	\label{eq: Q probability}
		\prod_{(u, v) \in \Lambda(t)} N_{u, v}(n) \bigg/ \sum_{t'} \prod_{(u, v) \in \Lambda(t')} N_{u, v}(n).
	\ee
	Now we take expectation on both sides of (\ref{eq: conditional law on transition counts}) and let $n\to\infty$ in (\ref{eq: Q probability}). This completes the proof of (\ref{eq: spanning tree probability}), since
	\[
		\lim_{n\to\infty}\frac{N_{u, v}(n)}{n} = \lim_{n\to\infty} \frac{N_{u}(n)}{n}
 \cdot \frac{N_{u, v}(n)}{N_{u}(n)}
		= \pi(u) p(u, v) = \frac{w(u, v)}{\sum_{u\in V} w(u)} \quad\P_x\text{-a.s.},
	\]
where $N_{u}(n):= \sum_{v \in V} N_{u, v}(n)$, and $\pi(u) := w(u) \big/ \sum_{u \in V} w(u)$ represents the stationary distribution of the random walk.
\end{proof}

A notable feature of our proof is that we never reversed time. Thus, essentially the same proof yields the following, where $\overleftarrow p$ denotes the probability transitions of the reversed Markov chain and $-e$ denotes the reversal of $e$:

\begin{theorem}[Directed Aldous--Broder]
\label{thm: Directed Aldous--Broder}
	Let $X$ be an irreducible Markov chain on a finite state space. Let $G$ be the corresponding directed graph. Let $X_0 = x$ and stop $X$ at the cover time. Let $T$ be the arborescence formed by all the first-entrance edges. Then there is a constant $c$ such that for every arborescence $t$ of $G$ rooted at $x$, $\P[T = t] = c \prod_{e \in t} \overleftarrow p(-e)$. 
\end{theorem}

\begin{proof}
It remains to note that $\prod_{(u, v) \in t} \pi(u) p(u, v) = 
\prod_{(u, v) \in t} \pi(v) \overleftarrow p(v, u) = \prod_{v \ne x} \pi(v) \cdot \prod_{e \in t} \overleftarrow p(-e)$.
\end{proof}

It turns out that one can also deduce the Markov chain tree theorem in a similar fashion. To this end, we introduce a random walk measure $\mu_n$ on the loop space. Let $X$ be an irreducible Markov chain on a finite state space with stationary probability measure $\pi$. Let $G$ be the corresponding directed graph. Recall that $\mathcal P_n^{x,y}$ denotes the set of all paths in $G$ from $x$ to $y$ with length $n$. Let $$\mcL_n:= \bigcup_{x \in V(G)} \mathcal P_n^{x,x}$$ be the space of all loops of length $n$ in $G$. If $r$ is the period of the Markov chain $X$, then $\mcL_n$ is nonempty iff $n \in r \mathbb N$. For simplicity in what follows, we will assume that $X$ is an aperiodic chain, i.e., $r=1$, although the proof is similar for $r>1$. Assuming $r=1$, we define a probability measure $\mu_n$ on $\mcL_n$ by assigning to each path $\gamma = (v_0, \ldots,v_n=v_0) \in \mcL_n$ the probability
	$$
	\mu_n(\gamma) := \frac{1}{Z_n} p(\gamma) = \frac{1}{Z_n} \prod_{k=0}^{n-1} p(v_k,v_{k+1}) ,
	$$
	where $Z_n$ is the normalizing constant. As $n$ goes to infinity, the ergodic theorem implies that $ Z_n \mu_n(\mathcal P_n^{x,x}) \to \pi(x)$. Thus, we have $Z_n \to 1$ and $\mu_n(\mathcal P_n^{x,x}) \to \pi(x)$.

\begin{corollary}[Markov Chain Tree Theorem]
\label{thm: MCTT}
	Let $X$ be an irreducible Markov chain on a finite state space with stationary probability measure $\pi$. Let $G$ be the corresponding directed graph. Then there is a constant $c$ such that for every $x$, $\pi(x) = c \sum_{\text{root}(t) = x} \prod_{e \in t} p(e)$, where the sum is over arborescences $t$ of $G$ directed \emph{toward} the root $x$. 
\end{corollary}

\begin{proof}
	The proof here follows closely that of Theorem \ref{thm: Aldous--Broder}. The main difference, however, is that we sample the path $(X_0,\ldots,X_n)$ not from the network random walk, but from the random walk measure $\mu_n$ on the loop space $\mcL_n$. Notice that the definitions of $N(n)$ and $D_n$ below (\ref{eq: transition counts}) still make sense as random variables on the loop space $\mcL_n$. In this case, the start of $D_n$ coincides with its end, and the indegree of each vertex in $D_n$ is equal to its outdegree. Let $D_n^\rmo$ be the same balanced multi-digraph as $D_n$, except that it does \emph{not} distinguish any start or end. All Eulerian paths in $D_n^\rmo$ are Eulerian circuits, which we consider as loops \emph{with} starts: different Eulerian paths may have different starts in $V(D_n^\rmo)$.
	
	 Now note that the law of random loops $\mu_n$, conditioned on the transition counts $N(n)$, has the same distribution as the law $\mathbb Q_{D_n^\rmo}$ of a uniformly random Eulerian path in $D_n^\rmo$. Also note that the probability under $\mathbb Q_{D_n^\rmo}$ for a certain arborescence of $D_n^\rmo$ to occur is no longer uniform, but proportional to the indegree of its root, according to (\ref{eq: number of Eulerian paths}). Therefore, for any arborescence $t$ of $G$ directed \emph{away} from its root, we have
	 $$
	 \mu_n[T_n=t \mid N(n)] = N_{\text{root}(t)}(n)\prod_{(u, v) \in t} N_{u, v}(n) \bigg/ \sum_{t'} N_{\text{root}(t')} \prod_{(u, v) \in t'} N_{u, v}(n),
	 $$
	 where $T_n$ is the arborescence formed by all the first-entrance edges.
	 By a similar calculation as in the proof of Theorem \ref{thm: Aldous--Broder} and \ref{thm: Directed Aldous--Broder}, we conclude that there exists a constant $c$ such that
		$$
		\lim_{n\to\infty} \mu_n[T_n=t] = c \prod_{e \in t} \overleftarrow p(-e).
		$$
		Thus
		$$
		\pi(x) = \lim_{n\to\infty} \mu_n(\mathcal P_n^{x,x}) = \lim_{n\to\infty} \sum_{\text{root}(t) = x} \mu_n[T_n=t] = c \sum_{\text{root}(t) = x} \prod_{e \in t} \overleftarrow p(-e).
		$$
		We complete the proof by reversing the process.
\end{proof}


Another consequence of Theorem \ref{thm: Directed Aldous--Broder} is the directed Wilson's algorithm. Due to the spatial Markov property of $\mathsf{UST}$ for digraphs, it suffices to show the following analogue of Corollary \ref{cor: path in UST has the same law as LERW}.  
\begin{corollary}\label{cor: path in UST has the same law as LERW for digraphs}
	Let $X$ be an irreducible Markov chain on a finite state space. Let $G$ be the corresponding directed graph. Let $x$ be a vertex in the digraph $G$ and  $\mathcal{T}_x$ be the collection of spanning trees of $G$ such that all edges on the trees are directed \emph{toward} $x$. For $t\in\mathcal{T}_x$, write $\Psi(t)=\prod_{e\in t}p(e)$. Let $\mathsf{UST}=\mathsf{UST}(G,x)$ be the uniform spanning tree measure, i.e., the probability measure on $\mathcal{T}_x$ such that  $\mathsf{UST}(\{t\})\propto \Psi(t)$. For a vertex $y\neq x$, and let $\gamma_{y,x}$ be the unique directed path in the $\mathsf{UST}$ from $y$ to $x$. Then $\gamma_{y,x}$ has the same law as the loop-erased random walk from $y$ to $x$. 
\end{corollary}
\begin{proof}
	Let $\overleftarrow{\mathbb{P}}_x$ denote the law of the reversed chain started at $x$. Let $r=(r_0,r_1,\dots,r_n)$ be a directed path in $G$ and let $-r=(r_n,r_{n-1},\dots,r_0)$ be the reversed path in the reversed graph. Similar to \eqref{eq:wts}, we let $\overleftarrow{p}(-r)=\prod_{j=1}^{n}\overleftarrow{p}(r_j,r_{j-1})$. 
	
		For a directed self-avoiding path $w$ from $y$ to $x$ in $G$, let $\mathscr{R}_{y,x}(w)$ be the set of directed paths $r=(r_0,r_1,\dots,r_n)$ on $G$ from $y$ to $x$ such that 
	\begin{itemize}
		\item $r_0=y,r_n=x$, and $\forall i\geq 1 \enspace r_i\neq y$,
		\item the loop-erasure of the path $r$ is $w$. 
	\end{itemize} 
Thus, $-\mathscr{R}_{y,x}(w):=\{-r\st  r\in \mathscr{R}_{y,x}(w) \}$ is just the set of all possible walks in the reversed chain from $x$ to the first hit of $y$ such that the path from $x$ to $y$ in the tree  formed by first-entrance edges of the walk is $-w$. 

By Theorem \ref{thm: Directed Aldous--Broder}, for a directed self-avoiding path $w=(w_0,w_1,\dots,w_k)$ from $y$ to $x$ in $G$, 
	\begin{eqnarray}\label{eq: 0.1}
\mathbb{P}[\gamma_{y,x}=w]&=& \sum_{r\in\mathscr{R}_{y,x}(w)} \overleftarrow{p}(-r)=\sum_{r\in\mathscr{R}_{y,x}(w)}  \frac{\pi(y)}{\pi(x)}p(r)\nonumber\\
&=& \frac{\pi(y)}{\pi(x)}\prod_{j=1}^{k}p(w_{j-1},w_j)\prod_{j=1}^{k}\mathscr{G}_{A_{j-1}}(w_j,w_j),
\end{eqnarray}
where $A_j:=\{w_0,\dots,w_j\}$ for $j=0,\dots,k$ and $\mathscr{G}_Z(a,a)$ is the expected number of visits to $a$ strictly before hitting $Z$ by a random walk on $G$ started from $a$. The last step of \eqref{eq: 0.1} can be seen by decomposing the path in a unique way as the concatenation of  edges $(w_{j-1},w_j)$ and cycles at $w_{1},\ldots,w_n$, with the cycle at $w_j$ avoiding $A_{j-1}$.

Let $l_{y\to x} $ denote the loop-erased random walk from $y$ to $x$. Using the same technique of decomposing the path (or using formula (12.2.2) of \cite{Lawler1999}), one has that
\be\label{eq: 0.2}
\mathbb{P}[l_{y\to x} =w]=\prod_{j=1}^{k}p(w_{j-1},w_j)\prod_{j=0}^{k-1}\mathscr{G}_{A_{j-1}\cup\{x\} }(w_j,w_j),
\ee
where $A_{-1}:=\emptyset$.

Comparing \eqref{eq: 0.1} with \eqref{eq: 0.2}, we see that it suffices to show the following identity: for any simple directed path $w=(w_0,\dots,w_k)$ from $y$ to $x$,
\be\label{eq: 0.3}
\frac{\pi(y)}{\pi(x)}\prod_{j=1}^{k}\mathscr{G}_{A_{j-1}}(w_j,w_j)=\prod_{j=0}^{k-1}\mathscr{G}_{A_{j-1}\cup\{x\} }(w_j,w_j).
\ee

By Cramer's rule, $\mathscr{G}_Z(a,a)=\frac{\det(I-P)[Z\cup \{a\}]}{\det(I-P)[Z]}$, where $P$ is the transition matrix and for a matrix $A$ indexed by $V(G)$ and $Z\subset V(G)$,   the matrix $A[Z]$ is the matrix obtained from $A$ by deleting all rows and columns indexed by an element of $Z$. Using this, we may simplify \eqref{eq: 0.3} to 
\be\label{eq: 0.4}
\frac{\pi(y)}{\det(I-P)[y]}=\frac{\pi(x)}{\det(I-P)[x]}.
\ee
Since \eqref{eq: 0.4} is easy to show by applying Cramer's rule to the system of equations determining the stationary probability measure $\pi$, we are done. 
\end{proof}

\begin{remark}
	In \cite{Lawler1999}, Lawler proved another very interesting identity and used it to prove Wilson's algorithm:
	\[
 \forall a,b\notin Z \quad	\mathscr{G}_Z(b,b)\mathscr{G}_{Z\cup \{b\}}(a,a)=\mathscr{G}_Z(a,a)\mathscr{G}_{Z\cup \{a\}}(b,b).
	\]
	From this identity one can easily deduce that the function in (12.2.3) on page 200 of \cite{Lawler1999} is  a symmetric function. The directed Wilson's algorithm then follows from this symmetry easily; see  (12.7.1) on page 212 of \cite{Lawler1999}. 
\end{remark}

\subsection{Infinite Networks}

Suppose $G$ is a locally finite,  infinite connected graph. An \notion{exhaustion} of $G$ is a sequence of finite connected subgraphs $G_n=(V_n,E_n)$ of $G$ such that $G_n\subset G_{n+1}$ and $G=\bigcup G_n$. Suppose that $V_n$ induces $G_n$, i.e., $G_n$ is the maximal subgraph of $G$ with vertex set $V_n$. Let $G_n^*$ be the graph formed from $G$ by contracting all vertices outside $V_n$ to a new vertex $\partial_n$. Let $T_n$ be a sample of $\mathsf{UST}(G_n^*)$. Then the \notion{wired uniform spanning forest} is the weak limit of $T_n$. If we orient $T_n$ toward $\partial_n$ and then take the weak limit, we get the \notion{oriented wired uniform spanning forest} of $G$. For details, see \cite{BLPS2001} or \cite[Chapter 10]{LP2016}.

Wilson's algorithm \cite{Wilson1996} is another efficient way of sampling uniform spanning tree  for finite connected graph using loop erasure of random walks. It can be applied to recurrent networks directly. For transient networks, Wilson's algorithm can also be applied with a simple modification. This is called Wilson's method rooted at infinity; see \cite{BLPS2001} for details.

 The Aldous--Broder algorithm can also be used to sample the wired uniform spanning forest on recurrent networks directly. However,
 the extension of the Aldous--Broder algorithm for sampling the wired uniform spanning forest on transient networks was found much later by 
Hutchcroft \cite{Hutchcroft2018interlacements} using the interlacement process.   

  Here we simply state a similar  interlacement using last-exit edges. It is not directly related to our reverse Aldous--Broder algorithm. Instead, it relates to the process created from the stationary random walk on the nonpositive integers that we recalled in our introduction. The interested reader can refer to  \cite{Hutchcroft2018interlacements} for details on the interlacement process and the interlacement Aldous--Broder algorithm.

 \begin{theorem}\label{thm: interlacement reverse AB}
 	Let $G$ be a transient, connected, locally finite network, let $\mathscr{I}$ be the interlacement process on $G$, and let $t\in\mathbb{R}$. For each vertex $v$ of $G$, let $\lambda_t(v)$ be the largest time before  $t$ such that there exists a trajectory $\bigl(W_{\lambda_t(v)},\lambda_t(v)\bigr)\in\mathscr{I}$ passing through $v$, and let $e_t(v)$ be the oriented edge of $G$ that is traversed by the trajectory $W_{\lambda_t(v)}$ as it leaves $v$ for the last time before $t$. Then
 	$
 	\{ e_t(v) \st v\in V \}
 	$
 	has the law of the oriented wired uniform spanning forest of $G$.
 \end{theorem}

 The proof of Theorem \ref{thm: interlacement reverse AB} is simply an analogue to that of \cite[Theorem 1.1]{Hutchcroft2018interlacements}. This version can be used in
 place of that used by Hutchcroft and is perhaps more natural.

 \newcommand\noopsort[1]{}
 \bibliography{revAB-reference}
 \bibliographystyle{plain}

\end{document}